\documentclass[nospthms,envcountsect]{svmult} 
\usepackage{makeidx}
\makeindex
\usepackage{latexsym}
\usepackage{amsfonts}
\usepackage{amsmath}
\usepackage{amssymb}
\usepackage{amscd}
\usepackage{multicol}    
\usepackage{enumerate}
    \newtheorem{theorem}     {Theorem} [section]
    \newtheorem{lemma}      [theorem]       {Lemma}
    \newtheorem{corollary}  [theorem]       {Corollary}
    \newtheorem{proposition}[theorem]       {Proposition}
    \newtheorem{definition} [theorem]       {Definition}
    \newtheorem{conjecture} [theorem]       {Conjecture}

\begin{document}
\catcode`@=11
\atdef@ I#1I#2I{\CD@check{I..I..I}{\llap{$\m@th\vcenter{\hbox
  {$\scriptstyle#1$}}$}
  \rlap{$\m@th\vcenter{\hbox{$\scriptstyle#2$}}$}&&}}
\atdef@ E#1E#2E{\ampersand@
  \ifCD@ \global\bigaw@\minCDarrowwidth \else \global\bigaw@\minaw@ \fi
  \setboxz@h{$\m@th\scriptstyle\;\;{#1}\;$}%
  \ifdim\wdz@>\bigaw@ \global\bigaw@\wdz@ \fi
  \@ifnotempty{#2}{\setbox@ne\hbox{$\m@th\scriptstyle\;\;{#2}\;$}%
    \ifdim\wd@ne>\bigaw@ \global\bigaw@\wd@ne \fi}%
  \ifCD@\enskip\fi
    \mathrel{\mathop{\hbox to\bigaw@{}}%
      \limits^{#1}\@ifnotempty{#2}{_{#2}}}%
  \ifCD@\enskip\fi \ampersand@}
\catcode`@=\active

\renewcommand{\labelenumi}{\alph{enumi})}
\newcommand{\chr}{\operatorname{char}}
\newcommand{\isom}{\stackrel{\sim}{\longrightarrow}}
\newcommand{\Aut}{\operatorname{Aut}}
\newcommand{\Shv}{\operatorname{Shv}}
\newcommand{\Hom}{\operatorname{Hom}}
\newcommand{\End}{\operatorname{End}}
\newcommand{\HOM}{\operatorname{{\mathcal H{\frak{om}}}}}
\newcommand{\Mod}{\operatorname{Mod}}
\newcommand{\EXT}{\operatorname{\mathcal E{\frak xt}}}
\newcommand{\Tot}{\operatorname{Tot}}
\newcommand{\Ext}{\operatorname{Ext}}
\newcommand{\Gal}{\operatorname{Gal}}
\newcommand{\cosk}{\operatorname{cosk}}
\newcommand{\Pic}{\operatorname{Pic}}
\newcommand{\Spec}{\operatorname{Spec}}
\newcommand{\trdeg}{\operatorname{trdeg}}
\newcommand{\holim}{\operatornamewithlimits{holim}}
\newcommand{\im}{\operatorname{im}}
\newcommand{\hyp}{{\rm hyp}}
\newcommand{\coim}{\operatorname{coim}}
\newcommand{\coker}{\operatorname{coker}}
\newcommand{\gr}{\operatorname{gr}}
\newcommand{\id}{\operatorname{id}}
\newcommand{\Br}{\operatorname{Br}}
\newcommand{\cd}{\operatorname{cd}}
\newcommand{\CH}{CH}
\newcommand{\Alb}{\operatorname{Alb}}
\renewcommand{\lim}{\operatornamewithlimits{lim}}
\newcommand{\colim}{\operatornamewithlimits{colim}}
\newcommand{\rk}{\operatorname{rank}}
\newcommand{\codim}{\operatorname{codim}}
\newcommand{\NS}{\operatorname{NS}}
\newcommand{\cone}{{\rm cone}}
\newcommand{\rank}{\operatorname{rank}}
\newcommand{\ord}{{\rm ord}}
\newcommand{\f}{{\cal F}}
\newcommand{\g}{{\cal G}}
\newcommand{\p}{{\cal P}}
\newcommand{\N}{{\mathbb N}}
\newcommand{\A}{{\mathbb A}}
\newcommand{\Z}{{{\mathbb Z}}}
\newcommand{\Q}{{{\mathbb Q}}}
\newcommand{\R}{{{\mathbb R}}}
\newcommand{\B}{{\mathbb Z}^c}
\renewcommand{\H}{{{\mathbb H}}}
\renewcommand{\P}{{{\mathbb P}}}
\newcommand{\F}{{{\mathbb F}}}
\newcommand{\m}{{\mathfrak m}}
\newcommand{\Sm}{{\text{\rm Sm}}}
\newcommand{\Sch}{{\text{\rm Sch}}}
\newcommand{\et}{{\text{\rm et}}}
\newcommand{\eh}{{\text{\rm eh}}}
\newcommand{\Zar}{{\text{\rm Zar}}}
\newcommand{\Nis}{{\text{\rm Nis}}}
\newcommand{\tr}{\operatorname{tr}}
\newcommand{\tor}{{\text{\rm tor}}}
\newcommand{\PreShv}{\text{\rm PreShv}}
\newcommand{\Div}{\operatorname{Div}}
\newcommand{\Ab}{{\text{\rm Ab}}}
\newcommand{\DD}{{\Z^c(0)}}
\renewcommand{\div}{\operatorname{div}}
\newcommand{\corank}{\operatorname{corank}}
\renewcommand{\O}{{\cal O}}
\newcommand{\C}{{\cal C}}
\renewcommand{\p}{{\mathfrak p}}
\newcommand{\proof}{\noindent{\it Proof. }}
\newcommand{\proofend}{\hfill $\Box$ \\}
\newcommand{\rem}{\noindent {\it Remark. }}
\newcommand{\example}{\noindent {\bf Example. }}
\newcommand{\ar}{{\text{\rm ar}}}
\newcommand{\del}{{\delta}}

\title*{Arithmetic homology and an integral version of Kato's conjecture}
\author{Thomas Geisser\thanks{Supported in part by NSF grant no.0556263}}
\institute{University of Southern California}

\maketitle

\begin{abstract}
We define an integral Borel-Moore homology theory over finite fields, called
arithmetic homology, and an integral version of Kato homology. Both types of
groups are expected to be finitely generated, and sit in a
long exact sequence with higher Chow groups of zero-cycles.
\end{abstract}

\section{Introduction}
For a separated scheme $X$ of finite type over a field $k$, let
$\Z^c(0)(X)$ be Bloch's complex of relative zero-cycles, generated
in degree $-i$ by cycles of dimension $i$ on $X\times \Delta^i$ in
good position. The higher Chow groups $CH_0(X,i)$ are defined as the
homology $H_i(\Z^c(0)(X))$ of $\Z^c(0)(X)$. Varying $X$, we obtain a
complex of etale sheaves $\Z^c(0)$ and can consider its etale
hypercohomology $H_i(X_\et,\Z^c(0))$. If $k$ is algebraically
closed, then it is a consequence of the Beilinson-Lichtenbaum
conjecture that $CH_0(X,i)\cong H_i(X_\et,\Z^c(0))$
\cite[Thm.3.1]{ichdual}, but in general, these groups are not
isomorphic. Over a finite field, Kato \cite{kato} defined for each
$m$ a complex with homology $H_i^K(X,\Z/m)$ and conjectured that if $X$ is smooth,
proper and connected, then $H_i^K(X,\Z/m)$ vanishes for $i>0$, and
$H_0^K(X,\Z/m)\cong\Z/m$. Jannsen and Saito \cite{jannsensaitoalt} observed
that Kato homology measures the difference
between the finite coefficient versions $CH_0(X,i,\Z/m)$ and
$H_i(X_\et,\Z^c/m(0))$, and proved Kato's conjecture assuming resolution of
singularities \cite{jannsensaito}.

In this paper, we construct Borel-Moore homology groups
$H_i^c(X_\ar,\Z)$, which are a substitute for the (pathological) etale
higher Chow groups $H_i(X_\et,\Z^c(0))$, and define an integral
version of Kato's complex whose homology groups $H_i^K(X,\Z)$ measure the
difference between $CH_0(X,i)$ and $H_i^c(X_\ar,\Z)$. The analog of
Kato's conjecture is that if $X$ is smooth, proper and connected,
then $H_i^K(X,\Z)$ vanishes for $i>0$, and $H_i^K(X,\Z)=\Z$.

The Borel-Moore homology theory, which we call arithmetic homology,
is constructed by applying the Weil-etale formalism of Lichtenbaum to
$\Z^c(0)$: Given a scheme $X$ over the finite field $\F_q$,
$H_i^c(X_\ar,\Z)$ is the $-i$th cohomology group of the
complex $R\Gamma_GR\Gamma(\bar X_\et,\Z^c(0))[1]$,
where $\bar X=X\times_{\F_q}\bar \F_q$, and $G$ is the Weil group of $\F_q$
acting on $\bar \F_q$.
These groups are expected to be finitely generated. More precisely, we
have

\begin{theorem}\label{introthm}
Let $X$ be smooth and proper. Assume resolution of singularities up to the
dimension of $X$ and the Beilinson-Lichtenbaum conjecture (see below).
Then the following statements are equivalent:
\begin{enumerate}
\item $CH_0(X,i)_\Q=0$ for all $i>0$.
\item There is an isomorphism $CH_0(X,i)\cong H_{i+1}^c(X_\ar,\Z)$ of
finitely generated abelian groups for all $i\geq 0$.
\item There are short exact sequences for all $i$,
$$ 0\to CH_0(\bar X,i+1)_G\to CH_0(X,i) \to CH_0(\bar X,i)^G\to 0.$$
\end{enumerate}
\end{theorem}

The vanishing of the Chow groups is a special case of Parshin's conjecture
$K_i(X)_\Q=0$ for $i>0$, and in particular it follows from Tate's conjecture
together
with Beilinson's conjecture that rational and homological equivalence
agree up to torsion \cite{ichtate} over finite fields,
or from finite dimensionality of
smooth and projective schemes over $\F_q$ in the sense of
Kimura-O'Sullivan \cite{ichparshin}.

The integral analog of Kato homology comes into play when
comparing higher Chow groups with arithmetic homology for not
necessarily smooth or proper schemes over $\F_q$: Consider the complex
$$\cdots\to\oplus_{x\in X_{(s)}}K_s^M( k(x)\otimes_{\F_q}\bar\F_q)_G
\to \oplus_{x\in X_{(1)}} (k(x)\otimes_{\F_q}\bar\F_q)^\times_G\to \cdots
\to \oplus_{x\in X_{(0)}}\Z, $$
where $X_{(s)}$ denote the points of $X$ of dimension $s$, and
$K_s^M( k(x)\otimes_{\F_q}\bar\F_q)_G$ is the group of Frobenius
coinvariants of the Milnor $K$-group of the finite product
of fields $ k(x)\otimes_{\F_q}\bar\F_q$.
The maps in the complex are induced by boundary maps in localization
sequences of higher Chow groups via the identification
$K_s^M( k(x)\otimes_{\F_q}\bar\F_q)
\cong H^s(k(x)\otimes_{\F_q}\bar\F_q,\Z(s))$.
The group $H_i^K(X,\Z)$ is the $i$th homology of this complex.
Assuming the Beilinson-Lichtenbaum conjecture, we have a long exact sequence
$$ \cdots \to H_i^K(X,\Z)\stackrel{\times m}{\longrightarrow}
H_i^K(X,\Z) \to H_i^K(X,\Z/m)\to \cdots, $$
which justifies calling the groups $H_i^K(X,\Z)$ an integral version
of Kato homology. As an analog of Kato's conjecture, we propose
the following

\begin{conjecture}\label{po1}
If $X$ is smooth, proper and connected, then
$H_i^K(X,\Z)=0$ for $i>0$ and $H_0^K(X,\Z)\cong\Z$.
\end{conjecture}

It is easy to use the known results on the (torsion) Kato conjecture
to show that the conjecture is true in degree $0$, and that
$H_i^K(X,\Z)\cong CH_0(X,i)_\Q$ for smooth and proper $X$ and $i=1,2$.
Regarding the conjecture
and the relationship between higher Chow groups and arithmetic
homology, we show

\begin{theorem}\label{intro2}
Assuming resolution of singularities and the Beilinson-Lichten\-baum
conjecture, the following statements are equivalent:

a) For every smooth and proper $X$ over $\F_q$, $CH_0(X,i)_\Q=0$ for $i>0$.

b) Conjecture \ref{po1} holds, and for every separated scheme of finite
type $X$ over $\F_q$, there is a long exact
sequence of finitely generated groups
$$\cdots \to CH_0(X,i)\to H_{i+1}^c(X_\ar,\Z)
\to H_{i+1}^K(X,\Z) \to CH_0(X,i-1) \to \cdots .$$
\end{theorem}

The exact sequence of the Theorem exists unconditionally in
degrees $i\leq 1$.
We show that Conjecture \ref{po1} holds for curves, i.e. there is
an exact sequence
\begin{equation}\label{AA}
0\to  (k(C)\otimes_{\F_q}\bar\F_q)^\times_G\to \bigoplus_{x\in C_{(0)}}
\Z\to \Z\to 0.
\end{equation}
for smooth and proper curves $C$.
Since the sequence comparing Weil-etale to etale cohomology \cite{ichweil}
gives a short exact sequence
$$ 0\to (k(C)\otimes_{\F_q}\bar\F_q)^\times_G\to k(C)^\times_\Q\to
\Br k(C)\to 0,$$
we obtain $ (k(C)\otimes_{\F_q}\bar\F_q)^\times_G\otimes\Q/\Z\cong \Br k(C)$,
and
\eqref{AA} is an integral version of the classical short exact sequence
for the Brauer group of $k(C)$.
For an arbitrary curve $C$ over $\F_q$, we obtain
$H_0^c(C_\ar,\Z)\cong H_0^K(C,\Z)$,
there is a short exact sequence
$$ 0\to CH_0(C)\to H_1^c(C_\ar,\Z)\to H_1^K(C,\Z)\to 0,$$
and $CH_0(C,i)\cong H_{i+1}^c(C_\ar,\Z)$ for $i\geq 1$.
The latter groups are finitely generated for $i=1$ and zero
for $i>1$. For curves, Kato homology is easy to calculate recursively;
for example if $C$ is proper with dual graph $\Gamma_C$, then
$H_i^K(C,\Z)=H_i(\Gamma_C,\Z)$.

Arithmetic homology can be applied to study abelian
class field theory of proper schemes.
The main observation is that the group $H_1^c(X_\ar,\Z)$
(which is conjecturally finitely generated) becomes isomorphic
to the abelianized fundamental group $\pi_1(X)^{ab}$ after profinite
completion, and the reciprocity map factors as
$$CH_0(X)\xrightarrow{\theta_X} H_1^c(X_\ar,\Z)\to \pi_1(X)^{ab}.$$

\medskip

{\it Notation:} For an abelian group $A$, $A^\wedge=\lim A/n$ is the
profinite completion, and $A^*=\Hom(A,\Q/\Z)$ the Pontrjagin dual.
All schemes over a field are separated and of finite type;
from section 3 we fix a finite field as the base field.

\medskip
{\it Acknowledgements:}  Parts of this paper were
written while the author visited the University of Tokyo,
and we thank T.Saito and the University for the inspiring atmosphere
they provided. We are indebted to S.Saito for inspiring discussions.

\section{Higher Chow groups}
For a scheme $X$ over a field $k$, Bloch's higher Chow complex
$z_n(X,*)$ is defined as follows \cite{bloch}. In degree $-i$,
it is  be the free abelian group $z_n(X,i)$ generated by cycles of
dimension $n+i$ on $X\times_k \Delta^i$ which meet all faces
properly. The differentials are given by taking the alternating sum of
intersection with face maps.
Higher Chow groups $CH_n(X,i)$ are defined as the homology of this
complex. We let ${\B}(n)_X=z_n(-,*)[2n]$ be the complex of etale
sheaves on $X$ with $z_n(-,-2n-i)$ in degree $i$.
For a proper map $f:X\to Y$ we have a push-forward
$f_*{\B}(n)_X\to {\B}(n)_Y$, and for a flat,
equidimensional map $f:X\to Y$ of relative dimension $d$, we have a
pull-back $f^*{\B}(n)_Y\to {\B}(n+d)_X[-2d]$.
If $X$ is smooth of dimension $d$, then there is a quasi-isomorphism of
complexes of Zariski sheaves ${\B}(n)\cong \Z(d-n)[2d]$, where the right hand
side is the motivic complex of Voevodsky \cite{voevodsky}.
For a finitely generated field $F$ over $k$, we define
$CH_n(F,i)=\colim_U CH_n(U,i)$, where
the colimit runs through $U$ of finite type over $k$ with field of
functions $F$.
If $F$ has transcendence degree $d$ over $k$, then
$CH_n(F,i)\cong H^{2d-2n-i}(F,\Z(d-n))$, where the right hand side
is motivic cohomology, and vanishes for $i<d-n$, and agrees with
$K_{d-n}^M(F)$ for $i=d-n$.
As a formal consequence of localization for higher Chow groups,
one obtains an isomorphism
\begin{equation}\label{indexing}
CH_n(X,i-2n)=H_i(X_\Zar,\Z^c(n))=:H_i^c(X_\Zar,\Z(n))
\end{equation}
and spectral sequences
\begin{equation}\label{niveau}
E^1_{s,t}=\oplus_{x\in X_{(s)}}H^{s-t}(k(x),\Z(s-n))
\Rightarrow H_{s+t}^c(X_\Zar,\Z(n)).
\end{equation}
In particular, $H_i^c(X_\Zar,\Z(n))=0$ for $i<n$.
For an abelian group $A$, we define
\begin{align}\label{define}
H_i^c(X_\et,A(n))=H^{-i}R\Gamma(X_\et,A\otimes \Z^c(n)).
\end{align}
If $A=\Q$, then $H_i^c(X_\et,\Q(n))=H_i^c(X_\Zar,\Q(n))=CH_n(X,i-2n)_\Q$.
In \cite{ichdual}, we proved that for every integer $m$ and every scheme
$f:X\to k$ over a perfect field $k$, there is a quasi-isomorphism
\begin{equation}\label{duali}
Rf^!\Z/m\cong \B/m(0).
\end{equation}
Here $Rf^!$ is the extraordinary inverse image of SGA 4 XVIII.
Hence the etale homology groups $H_i^c(X_\et,\Z/m(0))$ with coefficients
$\Z^c/m(0)$ agree with usual etale homology groups
$H_i(X_\et,\Z/m):=H^{-i}(X_\et,Rf^!\Z/m)$ of Laumon \cite{laumon}.
For a finitely generated field $F$ over $k$, we define
$H_i^c(F_\et,A(n))=\colim_U H_i^c(U_\et,A(n))$.

The Beilinson-Lichtenbaum conjecture over $k$ in homological weight $n$
says that if $X$ is a smooth scheme of dimension $d$ over $k$, then the
canonical map
$$CH_n(X,i-2n)=H_i^c(X_\Zar,\Z(n))\to  H_i^c(X_\et,\Z(n))$$
is an isomorphism for $i\geq n+d-1$.
By considering cohomological dimension, it then must be an isomorphism
for all $i$ if $k$ is algebraically closed and $n\leq 0$.
As a special case, the conjecture implies that
$H^{d-n+1}(F_\et,\Z(d-n))=0$ for fields $F$ of transcendence degree $d$
over $k$; this statement is often called "Hilbert's Theorem 90".
Over an algebraically closed field, Suslin \cite{suslinetale}
shows that the groups $CH_0(X,i,\Z/m)$ and $H_i^c(X_\et,\Z/m(0))$ are
isomorphic for $\chr k\not| m$, but it is not clear that the isomorphism
is induced by the canonical change of topology map.
In \cite{marcI}, it is shown that
the canonical map is an isomorphism if $m$ is a power of the characteristic.
In \cite{ichdual}, we use Suslin's theorem to show the following result.

\begin{proposition}\label{locok1}
Let $k$ be a perfect field and $n\leq 0$.

a) (Localization) For a closed embedding $i:Z\to X$, the canonical map
$\B_Z(n)\to Ri^!\B_X(n)$ is a quasi-isomorphism.

b) (Homotopy formula) If $p:X\times \A^m\to X$ is the projection, then we have
a quasi-isomorphism of complexes of etale sheaves
$Rp_*\B_{X\times \A^m}(n)\cong \B_X(n-m)$.

c) (Niveau spectral sequence) There is a spectral sequence
\begin{equation}\label{niveauet}
E^1_{s,t}=\oplus_{x\in X_{(s)}}H^{s-t}(k(x)_\et,\Z(s-n)) \Rightarrow
H_{s+t}^c(X_\et,\Z(n)).
\end{equation}
In particular, $H_i^c(X_\et,\Z(n))=0$ for $i<n$.
\end{proposition}

The last statement follows because $s+t<n$ implies that $s<n$ or $t<0$,
both of which imply that $E_{s,t}^1$ vanishes.
Note that by the homotopy formula, the Beilinson-Lichtenbaum conjecture
in homological weight $n\leq 0$ implies the Beilinson-Lichtenbaum
conjecture in all weights less than $n$.

\section{Arithmetic homology with compact support}
We fix a finite field $\F_q$ with Galois group $\hat G=\Gal(\bar
\F_q/\F_q)$ and let $G\subset \hat G$ be the Weil group of $\F_q$, i.e.
the subgroup of $\hat G$ generated by the Frobenius endomorphism $\varphi$.
Given a separated scheme of finite type $X$ over $\F_q$, let $\bar
X=X\times_{\F_q}\bar \F_q$. For an abelian group $A$,
we define arithmetic homology groups with compact support
and coefficients in $A$ as the homology groups of the complex
$R\Gamma_GR\Gamma(\bar X_\et,A\otimes{\B}(n))[1]$,
where $G$ acts on $\bar \F_q$,
so that
$$H_i^c(X_\ar,A(n))=H^{1-i}R\Gamma_GR\Gamma(\bar X_\et,A\otimes {\B}(n)).$$
The Leray spectral sequence for composition of functors degenerates into
short exact sequences
\begin{equation}\label{lerayses1}
0\to H_i^c(\bar X_\et,\Z(n))_G \to H_i^c(X_\ar,\Z(n)) \to
H_{i-1}^c(\bar X_\et,\Z(n))^G\to 0 .
\end{equation}
In particular, $H_i^c(X_\ar,\Z(n))=0$ for $i<n\leq 0$
by Proposition \ref{locok1}.
If $X$ is smooth of dimension $d$, then the quasi-isomorphism
$\Z^c(n)\cong \Z(d-n)[2d]$ implies
\begin{equation}\label{smooth}
H_i^c(X_\ar,\Z(n))\cong H^{2d+1-i}(X_W,\Z(d-n)),
\end{equation}
where the right hand side are the Weil-etale cohomology groups
of \cite{licht, ichweil}.
(For arbitrary $X$ of finite type over $\F_q$, one has to replace etale
cohomology by eh-cohomology and Weil-etale cohomology by arithmetic
cohomology of \cite{ichweilII} to obtain well-behaved cohomology groups).
Recall from \cite[Thm.7.1]{ichweil} that for every smooth $X$ over $\F_q$,
there is a long exact sequence
\begin{equation}\label{alt}
\cdots\to H^i(X_\et,\Z(p))\to H^i(X_W,\Z(p))\to H^{i-1}(X_\et,\Q(p))\to\cdots.
\end{equation}
The same proof shows

\begin{theorem}\label{weilI}
There are long exact sequences
\begin{multline*}
\cdots\to H_{i-1}^c(X_\et,\Z(n)) \to H_i^c(X_\ar,\Z(n))\\
\to CH_n(X,i-2n)_\Q \to H_{i-2}^c(X_\et,\Z(n)) \to\cdots,
\end{multline*}
hence for every integer $m\geq 1$ an isomorphism
$$H_i^c(X_\et,\Z/m(n)) \cong  H_{i+1}^c(X_\ar,\Z/m(n)).$$
With rational coefficients, the long exact sequence splits, and
$$H_i^c(X_\ar,\Q(n)) \cong
\CH_n(X,i-2n)_\Q\oplus \CH_n(X,i-2n-1)_\Q.$$
\end{theorem}
\proofend

From now on, we assume that $n\leq 0$. Then
for a closed subscheme $Z$ of $X$ with open complement $U$ we have
by Proposition \ref{locok1} localization sequences
\begin{equation}\label{arithloc}
\cdots \to H_i^c(Z_\ar,\Z(n))\to
H_i^c(X_\ar,\Z(n)) \to H_i^c(U_\ar,\Z(n)) \to\cdots
\end{equation}

\begin{corollary} Let $n\leq 0$.

a) (Projective bundle formula) Let ${\mathbb P}^r_X\to X$ be a
projective bundle of relative dimension $r$. Then
$$H_i^c({\mathbb P}^r\times X_\ar,\Z(n))\cong
\oplus_{j=0}^r H_{i-2j}^c(X_\ar,\Z(n-j)).$$

b) (Homotopy formula) For every $r\geq 0$,
$$ H_i^c({\mathbb A}^r\times X_\ar,\Z(n))\cong H_{i-2r}^c(X_\ar,\Z(n-r)).$$
\end{corollary}

\proof
The homotopy formula follows from Prop. \ref{locok1} via \eqref{lerayses1},
and the projective bundle formula can be derived from this using
localization and induction.
\proofend

The Beilinson-Lichtenbaum conjecture over $\bar \F_q$ implies that,
for $n\leq 0$,
$R\Gamma(X_\ar,A\otimes {\B}(n))$ has the explicit representative
$${\B}(n)({\bar X})\otimes A[G] \xrightarrow{1-\varphi}
{\B}(n)({\bar X})\otimes A[G],$$
with $\varphi$ acting diagonally,
and that the sequence \eqref{lerayses1} takes the form
\begin{equation}\label{lerayses}
0\to CH_n(\bar X,i-2n)_G \to H_i^c(X_\ar,\Z(n)) \to
CH_n(\bar X,i-1-2n)^G\to 0.
\end{equation}
For example, $H_0^c(X_\ar,\Z(0))=CH_0(\bar X)_G$.

In the following, we restrict ourselves to the case $n=0$ and
write $H_i^c(X_\ar,A)$ for $H_i^c(X_\ar,A(0))$.
This is no loss of generality in view of the homotopy formula.

\example Let $P$ be a connected, zero-dimensional scheme over
$\F_q$. Then $CH_0(\bar P,i)=\Z^c$ for $i=0$ and
zero otherwise, where $c$ is the number of connected components
of $\bar P$. Since the Galois group permutes the connected
components of $\bar P$, \eqref{lerayses} implies
$H_i^c(P_\ar,\Z)\cong \Z$ for $i=0,1$, and zero otherwise.

\subsection{Arithmetic cohomology of fields}
For a finitely generated field $k$ over $\F_q$,
$k\otimes_{\F_q}\bar\F_q$ is a finite product of algebraic
extensions of $k$, finitely generated over $\bar \F_q$, on which $G$ acts.
We define arithmetic cohomology $H^i(k_\ar,A(n))$ as the cohomology of
$R\Gamma_GR\Gamma((k\otimes_{\F_q}\bar\F_q)_\et,A\otimes\Z(n))$. There exist
sequences \eqref{alt} comparing arithmetic cohomology to etale cohomology.
If the transcendence degree of $k$ over
$\F_q$ is $d$, then we can write $k$ as the colimit
of smooth schemes of dimension $d$ over $\F_q$, $k=\colim U$, and
$H^i(k_\ar,\Z(d))\cong \colim H^i(U_W,\Z(d))\cong
\colim H_{2d+1-i}^c(U_\ar,\Z(0))$ by \eqref{smooth}.
As a formal consequence of localization \eqref{arithloc},
we obtain

\begin{corollary}
For $n\leq 0$, there are spectral sequences
\begin{equation}\label{niveau1}
E^1_{s,t}=\oplus_{x\in X_{(s)}}H^{s-t}(k(x)_\ar,\Z(s-n))
\Rightarrow H_{s+t+1}^c(X_\ar,\Z(n)).
\end{equation}
\end{corollary}

\begin{lemma}\label{aroffields}
Assuming the Beilinson-Lichtenbaum conjecture in homological weight $d-s$,
then for every field $k$ with $d=\trdeg_{\F_q}k$ we have
isomorphisms
\begin{align*}
H^{s+1}(k_\ar,\Z(s))&\cong K_s^M(k\otimes_{\F_q}\bar\F_q)_G;\\
H^{s+2}(k_\ar,\Z(s))&\cong H^{s+1}( k\otimes_{\F_q}\bar\F_q,\Q/\Z(s))^G.
\end{align*}
If $s\geq d$, then the latter group vanishes, and we obtain
an exact sequence
\begin{equation}\label{3.6}
0\to K_s^M( k\otimes_{\F_q}\bar\F_q)_G \to K_s^M(k)_\Q\xrightarrow{\del}
H^{s+1}(k_\et,\Q/\Z(s)) \to 0,
\end{equation}
where $\del$ is the composition
$$ K_s^M(k)_\Q \to K_s^M(k)\otimes \Q/\Z\xrightarrow{sym}
H^s(k_\et,\Q/\Z(s))\xrightarrow{\cup e} H^{s+1}(k_\et,\Q/\Z(s)),$$
with $sym$ the symbol map, and
$e\in H^1((\F_q)_\et,\hat \Z)=\Hom_{cont}(\hat G,\hat \Z)$
the isomorphism sending the Frobenius endomorphism to $1$.
\end{lemma}

\proof
The first statements follow from the Leray spectral sequence
$$0\to H^{i-1}((k\otimes_{\F_q}\bar\F_q)_\et,\Z(s))_G
\to H^i(k_\ar,\Z(s))\to H^i((k\otimes_{\F_q}\bar\F_q)_\et,\Z(s))^G\to 0$$
and Hilbert's Theorem 90.
The exact sequence is \eqref{alt}, using the identifications
$H^s(k_\et,\Q(s))=K_s^M(k)_\Q$ and
$H^{s+2}(k_\et,\Z(s))\cong H^{s+1}(k_\et,\Q/\Z(s))$.
The description of $\del$ in \eqref{alt} is given in
in \cite[Thm.7.1a)]{ichweil}.
\proofend

The argument of the Lemma for $s=0,1$ gives, for any $k$, exact sequences
$$
0\to \Z \to\Q \to \Hom(G_k,\Q/\Z)\to
\Hom_G(\Gal(k\otimes_{\F_q}\bar\F_q),\Q/\Z) \to 0$$
\begin{equation}\label{ikmn}
0\to ( k\otimes_{\F_q}\bar\F_q)^\times_G \to (k^\times)_\Q \to \Br(k) \to
\Br( k\otimes_{\F_q}\bar\F_q)^G  \to 0.
\end{equation}

\begin{lemma}\label{iuy}
Under the hypothesis of the previous Lemma, there are short exact sequences
\begin{multline*}
0\to H^{s+1}(k_\ar,\Z(s))\to H^s(k,\Q(s))
\to H^{s+1}(k_\et,\Q/\Z(s))\to 0;
\end{multline*}
and
\begin{multline*}
0\to H^{s+1}(k_\ar,\Z(s))\stackrel{\times m}{\longrightarrow}
H^{s+1}(k_\ar,\Z(s))
\to H^{s+1}(k_\et,\Z/m(s))\to 0.
\end{multline*}
\end{lemma}

\proof
The first short exact sequence is \eqref{3.6}.
The long exact coefficient sequence for arithmetic
cohomology gives the second short exact sequence, because
the boundary map factors
$$\begin{CD}
H^s(k_\et,\Z/m(s))@>\delta >> H^{s+1}(k_\et,\Z(s))=0\\
@| @VVV \\
H^s(k_\ar,\Z/m(s))@>\delta >> H^{s+1}(k_\ar,\Z(s)),
\end{CD}$$
the upper right group is trivial by Hilbert's theorem 90.
\proofend

\section{Finite generation}
We fix a finite field $\F_q$.
\medskip

\noindent{\bf Conjecture $P_0(X)$:}
{\it For the smooth and proper scheme $X$
over $\F_q$, the group $CH_0(X,i)$ is torsion for all $i>0$.}

\medskip

This is a special case of Parshin's conjecture, because if $X$ has
dimension $d$, then
$$CH_0(X,i)_\Q=H^{2d-i}(X,\Q(d))=K_i(X)_\Q^{(d)}.$$
If Tate's conjecture holds and rational equivalence and homological
equivalence agree up to torsion for all $X$,
then Conjecture $P_0(X)$ holds for all $X$ \cite{ichtate}.

\begin{proposition}\label{4.0}
For a smooth and proper scheme $X$, the following statements are
equivalent:
\begin{enumerate}
\item Conjecture $P_0(X)$.
\item The groups $H_i^c(X_\ar,\Z)$ are finitely generated for
all $i$.
\item There are isomorphisms
$H_i^c(X_\ar,\Z)\otimes \Z_l\xrightarrow{\sim} \lim H_{i-1}(X_\et,\Z/l^r)$
for all $i$, and all $l$.
\item There are isomorphisms
$H_i^c(X_\ar,\Z)\otimes \Z_l\xrightarrow{\sim} \lim H_{i-1}(X_\et,\Z/l^r)$
for all $i$, and some $l$.
\end{enumerate}
\end{proposition}

\proof
a) $\Rightarrow$ b): If $X$ is smooth and proper, then
$H_i^c(X_\ar,\Z)\cong H^{2d+1-i}(X_W,\Z(d))$ by \eqref{smooth}.
By hypothesis and \eqref{alt}, for $j\leq 2d$,
$$H^j(X_W,\Z(d))\cong H^j(X_\et,\Z(d))\cong H^{j-1}(X_\et,\Q/\Z(d)),$$
which is finite for $j< 2d$ by the Weil-conjectures and Gabber \cite{gabber}.
By \cite{katosaito}, the group $H^{2d}(X_W,\Z(d))\cong H^{2d}(X_\et,\Z(d))$
agrees with $CH_0(X)$, and is finitely generated. Finally,
$H^{2d+1}(X_W,\Z(d))\cong \Z$ and $H^i(X_W,\Z(d))=0$ for $i>2d+1$
by \cite[Thm.7.5]{ichweil}.

b) $\Rightarrow$ c): Finite generation implies
$H_i^c(X_\ar,\Z)\otimes \Z_l \cong \lim H_i^c(X_\ar,\Z/l^r)\cong
\lim H_{i-1}(X_\et,\Z/l^r)$ by Theorem \ref{weilI} and \eqref{duali}.

d) $\Rightarrow$ a): The group
$\lim H_{i-1}(X_\et,\Z/l^r)=H^{2d-i}(X_\et,\Z_l(d))$ is
torsion for $i>1$ by the Weil-conjectures, and a) follows via
Theorem \ref{weilI}.
\proofend

\begin{proposition}\label{4.1}
The following statements are equivalent:
\begin{enumerate}
\item Conjecture $P_0(X)$ for all smooth and proper $X$.
\item The groups $H_i^c(X_\ar,\Z)$ are finitely generated for
all $i$ and all $X$.
\item There are isomorphisms
$H_i^c(X_\ar,\Z)\otimes \Z_l\xrightarrow{\sim} \lim H_{i-1}(X_\et,\Z/l^r)$
for all $i$, all $X$, and all $l$.
\item There are isomorphisms
$H_i^c(X_\ar,\Z)\otimes \Z_l\xrightarrow{\sim} \lim H_{i-1}(X_\et,\Z/l^r)$
for all $i$, all $X$, and some $l$.
\end{enumerate}
\end{proposition}

\proof
Following the proof of Prop. \ref{4.0}, it suffices to show that
finite generation
of $H_i^c(X_\ar,\Z)$ for smooth and proper $X$ implies finite generation
for all $X$. By localization, induction on the number of irreducible
components, and induction on the dimension, we reduce to the case
that $X$ is irreducible, and see that finite generation for $X$ and any
of its open subschemes $V$ are equivalent. Let $X'\to X$ be an
alteration such that $X'$ is an open subscheme of a smooth and proper
scheme $T$.
We can assume that there is an open subset $U\subseteq X$ such that
$U'=U\times_XX'\to U$ is Galois with group $A$.
Since $H_i^c(T_\ar,\Z)$ is finitely generated, so is $H_i^c(U'_\ar,\Z)$.
But there is a spectral sequence
$$E^2_{s,t}=H_s(A,H_t^c(U'_\ar,\Z(n)))\Rightarrow H_{s+t}^c(U_\ar,\Z(n)),$$
and group homology of a finite group with finitely generated
coefficients is finite. So $H_i^c(U_\ar,\Z)$, and hence
$H_i^c(X_\ar,\Z)$ is finitely generated.
\proofend

\begin{proposition}\label{degreenull}
Let $X$ be connected. If $X$ is proper, then the degree map induces an
isomorphism $H_0^c(X_\ar,\Z)\cong \Z$, and if $X$ is not proper,
then $H_0^c(X_\ar,\Z)=0$.
\end{proposition}

\proof
By \eqref{smooth} and \cite[Thm.7.5]{ichweil}, we have
$H_0^c(X_\ar,\Z)\cong H^{2d+1}(X_W,\Z(d))\cong \Z$ for
$X$ smooth and proper of dimension $d$.
In the general case, we can assume that $X$ is irreducible
by induction on the number of irreducible components.
By induction on the dimension, we obtain $H_0^c(U'_\ar,\Z)=0$
for the complement $U'$ of a non-empty closed subscheme in
a smooth and proper scheme $X'$.
With the argument of Prop. \ref{4.1} this implies that there is
an open subscheme $U$ of $X$ with connected complement
such that $H_0^c(U_\ar,\Z)=0$.
But then by induction on the dimension,
we get that $\Z=H_0^c((X-U)_\ar,\Z)\twoheadrightarrow H_0^c(X_\ar,\Z)
\to H_0^c((\F_q)_\ar,\Z)\cong \Z$ is an isomorphism.
\proofend

Consider the complex
\begin{equation}\label{rational}
\cdots\to \oplus_{x\in X_{(s)}}H^{s}(k(x),\Q(s))\to
\cdots \to \oplus_{x\in X_{(0)}}H^0(k(x),\Q(0)),
\end{equation}
arising as the $E_1$-terms and $d_1$-differentials in the
niveau spectral sequence \eqref{niveau}.
Let $\tilde H_i(X,\Q(0))=E^2_{i,0}(X)$
be the homology of this complex. Then in \cite{ichparshin}, we used
a theorem of Jannsen \cite{jannsen} to show that Conjecture $P_0(X)$ has two
independent faces:

\begin{proposition}\label{eqeqe}
Assume resolution of singularities. Then conjecture
$P_0(X)$ for all smooth and proper $X$ is equivalent to the following
statements for all smooth and proper $X$:
\begin{enumerate}
\item $CH_0(X,i)_\Q=0$ for $i>\dim X$.
\item $\tilde H_i(X,\Q(0))$ for all $i>0$ and $\tilde H_0(X,\Q(0))=\Q$.
\end{enumerate}
\end{proposition}

The first statement is equivalent to $H^i(k,\Q(d))=0$
for all fields $k$ of transcendence degree $d$ over $\F_q$ and
all $i\not=d$, or to $CH_0(X,i)_\Q\cong \tilde H_i(X,\Q(0))$
for all $X$ and all $i$. The second statement is equivalent to
the vanishing of $\tilde H_d(X,\Q(0))=\tilde H_{d-1}(X,\Q(0))=0$
for $d=\dim X>1$ \cite{ichparshin}.

\subsection{Special values of zeta-functions}
We can use Borel-Moore arithmetic homology to give formulas for
special values of zeta-functions at non-positive integers $n$.
There are other versions in Kahn \cite[Thm. 72]{kahnhandbook}
(up to powers of $p$),
and \cite[Conj. 1.4]{ichweilII} (assuming resolution of singularities).
Let $\chi(H_*^c(X_\ar,\Z(n)),e)$ be the Euler characteristic of the complex
\begin{equation}\label{euler}
\cdots\to H_{i+1}^c(X_\ar,\Z(n))\xrightarrow{\cap e} H_i^c(X_\ar,\Z(n))
\xrightarrow{\cap e} H_{i-1}^c(X_\ar,\Z(n))\to\cdots.
\end{equation}
The cap product is induced by the group cohomology product
with a generator $e$ of $H^1(G,\Z)\cong \Z$.

\begin{theorem}\label{fingen}
Under Conjecture $P_0(X)$ for all smooth and projective $X$,
the following statements hold for all $n\leq 0$ and all $X$ of
finite type over $\F_q$.
\begin{enumerate}
\item The groups $H_i^c(X_\ar,\Z(n))$ are finitely generated, and
$$\ord_{s=n}\zeta(X,s)=\sum_{i}i(-1)^i\cdot \rk H_i^c(X_\ar,\Z(n))=:\rho_n.$$
\item For $s\to n$,
$$ \zeta(X,s)\sim (1-q^{n-s})^{\rho_n}\cdot \chi(H_*^c(X_\ar,\Z(n)),e)^{-1}.$$
\end{enumerate}
\end{theorem}

\proof
By the homotopy formula, we can assume that $n=0$.
Using Prop. \ref{4.1} and the argument of
\cite[Thm.7.1]{ichweilII} we can reduce to the case where $X$ is smooth,
projective, and connected.
In this case $\ord_{s=0}\zeta(X,s)=\rho_0=-1$, and by Prop.\ref{4.0}b),
$$H_i^c(X_\ar,\Z)\otimes \hat\Z \cong H^{2d+1-i}(X_W,\Z(d))\otimes\hat\Z\\
\cong H^{2d+1-i}(X_\et,\hat\Z(d))$$
is the usual etale cohomology. By Poincare-duality,
$$H^{2d+1-i}(X_\et,\hat\Z(d))\cong \lim_m (H^i(X_\et,\Z/m)^*)
\cong H^i(X_\et,\Q/\Z)^*,$$
and since $D^*$ is torsion free for every divisible group $D$, we obtain
$$H^{2d+1-i}(X_\et,\hat\Z(d))_\tor \cong
\big((H^i(X_\et,\Q/\Z)/\Div)^*\big)_\tor
\cong (H^{i+1}(X_\et,\hat\Z)_\tor)^*.$$
It is easy to see that the map
$H^{2d}(X_\et,\hat\Z(d))/\tor\xrightarrow{\cup e}
H^{2d+1}(X_\et,\hat\Z(d))/\tor$
is an isomorphism. Since all other groups in \eqref{euler} are
finite, we obtain
\begin{multline*}
\chi(H_*^c(X_\ar,\Z),e)^{-1}
=\prod_i |H_i^c(X_\ar,\Z)_\tor|^{(-1)^{i-1}}
= \prod_i |H^{i+1}(X_\et,\hat\Z)_\tor|^{(-1)^{i-1}}\\
= \prod_i |H^i(X_W,\Z)_\tor|^{(-1)^i}
= \chi(H^*(X_W,\Z),e).
\end{multline*}
The latter is equal to the limit of $\zeta(X,s)(1-q^{-s})^{-1}$ for
$s\to 0$ by \cite[Thm.9.1]{ichweil}.
\proofend

\section{An integral version of Kato's conjecture}
Throughout this section we assume the validity
of the Beilinson-Lichtenbaum conjecture in weight $0$.
By Jannsen-Saito-Sato \cite[Thms. 1.6, 2.21]{JSS}, the complex
of $E_1^{s,-1}$-terms and differentials in \eqref{niveauet}
\begin{equation}\label{katoseq}
\cdots\to \oplus_{x\in X_{(s)}}H^{s+1}(k(x)_\et,\Z/m(s))\to
\cdots\to \oplus_{x\in X_{(0)}}H^{1}(k(x)_\et,\Z/m(0))\to 0
\end{equation}
is isomorphic to the complex defined by Kato in \cite{kato}.

\begin{conjecture}\label{katororign} (Kato)
For a smooth, proper, and connected scheme $X$
over a finite field, the homology groups $H_i^K(X,\Z/m)$ of \eqref{katoseq}
satisfy
$$H_i^K(X,\Z/m)=
\begin{cases}
\Z/m & i=0; \\
0  &i \not=0.
\end{cases}$$
\end{conjecture}

Comparing the niveau spectral sequence \eqref{niveau} and \eqref{niveauet},
one sees that there is a long exact sequence
\begin{equation}\label{jsseq}
\cdots\to CH_0(X,i,\Z/m)\to H_i(X_\et,\Z/m)\to
H_{i+1}^K(X,\Z/m) \to \cdots.
\end{equation}
Indeed, with finite coefficients,
the terms $E_1^{s,t}$ in \eqref{niveau} and \eqref{niveauet} agree for
$t\geq 0$ by the Beilinson-Lichtenbaum conjecture,
the term $E_1^{s,t}$ in \eqref{niveau} vanish for
$t<0$ by definition, the terms $E_1^{s,t}$ in \eqref{niveauet} vanish for
$t<-1$ by considering cohomological dimension, and the remaining terms
$E_1^{s,-1}$ in \eqref{niveauet} are given by the complex \eqref{katoseq}.

\subsection{The theorem of Jannsen and Saito}
\begin{theorem}\label{ssf}
(Jannsen, Saito \cite{jannsen, jannsensaito})
If $X$ is smooth, proper and connected, and
resolution of singularities holds for schemes up to dimension
$\min\{i,\dim X\}$, then $H_0^K(X,\Q/\Z)\cong \Q/\Z$,
and $H_i^K(X,\Q/\Z)=0$ for $i>0$.
\end{theorem}

The case $1\leq i\leq 3$ was previously treated by
Colliot-Th\'el\`ene \cite{colliet} and Suwa \cite{suwa}.
We can now prove Theorem \ref{introthm} of the introduction:

\proof (Thm.\ref{introthm})
a) $\Rightarrow$ b): By $P_0(X)$ and Thm. \ref{weilI}, the isomorphism
holds rationally. With torsion coefficients, it follows from
Thms. \ref{weilI}, \ref{ssf} and \eqref{jsseq}. Hence we get the isomorphism
by comparing coefficient sequences.
The groups are finitely generated by Prop. \ref{4.0}.

\noindent b) $\Rightarrow$ c) is \eqref{lerayses}, and
c) $\Rightarrow$ a) because
$CH_0(\bar X,i)_G\otimes \Q\cong CH_0(X,i)_\Q\cong CH_0(\bar X,i)^G\otimes\Q$
by a transfer argument.
\proofend

\begin{proposition}\label{5.3}
Assuming resolution of singularities, the following statements are equivalent:
\begin{enumerate}
\item $P_0(X)$ for all smooth and proper $X$.
\item $CH_0(X,i)$ is finitely generated for all $X$ and all $i$.
\item The canonical map
$CH_0(X,i)\otimes \Z_l\to \lim CH_0(X,i,\Z/l^r)$ is an isomorphism
for all $i$, all $l$, and all $X$.
\end{enumerate}
\end{proposition}

\proof
a) $\Rightarrow$ b): By the standard argument using localization,
it suffices to
show finite generation for smooth and proper $X$. In this case,
for $i>0$, by $P_0(X)$ and Theorem \ref{ssf},
\begin{multline*}
CH_0(X,i)\cong {}_\tor CH_0(X,i)\twoheadleftarrow CH_0(X,i+1,\Q/\Z)\\
\xrightarrow{\sim} H_{i+1}(X_\et,\Q/\Z)=H^{2d-i-1}(X_\et,\Q/\Z(d))
\end{multline*}
and this is finite by \cite{gabber} and the Weil-conjectures.
Finiteness of $CH_0(X)$ is a result of class field theory \cite{katosaito}.

\noindent b) $\Rightarrow$ c) is clear, and c) $\Rightarrow$ a) follows from
\begin{multline*}
CH_0(X,i)\otimes\Z_l \xrightarrow{\sim} \lim CH_0(X,i,\Z/l^r) \\
\xrightarrow{\sim} \lim H_i(X_\et,\Z/l^r)=\lim H^{2d-i}(X_\et,\Z/l^r(d)),
\end{multline*}
because the right hand group is torsion for smooth and proper $X$
and for $i>0$ by the Weil-conjectures.
\proofend

\subsection{The complex}
\begin{definition}
For a separated scheme of finite type $X$ over $\F_q$, we define
integral Kato-homology
$H_i^K(X,\Z)$ to be the homology of the complex $C^K(X)$ of
$E_1^{s,-1}$-terms and differentials in \eqref{niveau1}
\begin{equation}\label{gkatos}
\cdots\to \oplus_{x\in X_{(s)}}H^{s+1}(k(x)_\ar,\Z(s))
\to\cdots\to \oplus_{x\in X_{(0)}}H^1(k(x)_\ar,\Z(0))
\to 0.
\end{equation}
\end{definition}

By Lemma \ref{aroffields}, the complex \eqref{gkatos} can be rewritten as
\begin{equation}\label{tgfr}
\to \oplus_{x\in X_{(s)}}K_s^M( k(x)\otimes_{\F_q}\bar\F_q)_G
\to \cdots\to\oplus_{x\in X_{(1)}} (k(x)\otimes_{\F_q}\bar\F_q)^\times_G
\to \oplus_{x\in X_{(0)}}\Z .
\end{equation}

An inspection of the spectral sequence \eqref{niveau1} shows
that there is a map of complexes $C^K(X)\to H_0^c(X,\Z)$.
In particular, if $X$ is proper and connected, then we obtain from
Proposition \ref{degreenull} an augmentation $C^K(X)\to \Z$.
A closed embedding $Z\to X$ with open complement $U$
gives rise to a short exact sequence of complexes \eqref{tgfr},
hence a long exact sequence
$$ \cdots \to H_i^K(Z,\Z)\to H_i^K(X,\Z)\to H_i^K(U,\Z) \to \cdots.$$
The connection to Kato homology with finite coefficients is given as follows.

\begin{proposition}
There are long exact sequences
\begin{equation}\label{coecoe}
\cdots \to H_i^K(X,\Z)\stackrel{\times m}{\to}
H_i^K(X,\Z) \to H_i^K(X,\Z/m)\to \cdots ,
\end{equation}
and
\begin{equation}\label{relate}
\cdots \to H_i^K(X,\Z)\to \tilde H_i(X,\Q(0)) \to H_i^K(X,\Q/\Z)\to \cdots.
\end{equation}
\end{proposition}

\proof
By Lemma \ref{iuy}, there are short exact sequence of complexes
\begin{multline*}
0\to \bigoplus_{x\in X_{(*)}}H^{*+1}(k(x)_\ar,\Z(*))
\stackrel{\times m}{\longrightarrow}
\bigoplus_{x\in X_{(*)}}H^{*+1}(k(x)_\ar,\Z(*)) \\
\to \bigoplus_{x\in X_{(*)}}H^{*+1}(k(x)_\et,\Z/m(*))\to 0
\end{multline*}
and
\begin{multline*}
0\to \bigoplus_{x\in X_{(*)}}H^{*+1}(k(x)_\ar,\Z(*))
\to \bigoplus_{x\in X_{(*)}}H^*(k(x),\Q(*)) \\
\to \bigoplus_{x\in X_{(*)}}H^{*+1}(k(x)_\et,\Q/\Z(*))\to 0.
\end{multline*}
\proofend

\begin{corollary}\label{known}
If $X$ is smooth, proper and connected, then $H_0^K(X,\Z)\cong \Z$
and $H_i^K(X,\Z)\cong CH_0(X,i)_\Q$ for $i=1,2$.
\end{corollary}

\proof
Since $H^i(k(x),\Q(1))=0$ for $i\not=1$ and $H^i(k(x),\Q(0))=0$ for $i\not=0$,
an inspection of the niveau spectral sequence \eqref{niveau} shows that
$CH_0(X,i)_\Q\cong \tilde H_i(X,\Q(0))$ for $i\leq 2$
and for all $X$. The result follows with
the sequence \eqref{relate} and Theorem \ref{ssf}.
\proofend

\begin{theorem}\label{maincft}
a) For every $X$, there is a natural exact sequence
\begin{multline*}
CH_0(X,1)\to
H_2^c(X_\ar,\Z)\to H_2^K(X,\Z)\stackrel{d_0}{\longrightarrow}CH_0(X)\\
\to H_1^c(X_\ar,\Z)\to H_1^K(X,\Z)\to 0,
\end{multline*}
and the image of the map $d_0$ is torsion.

b) The statement in Proposition \ref{eqeqe}a) holds if and only
if there exists a long exact sequence
$$\cdots \xrightarrow{d_i} CH_0(X,i)\to H_{i+1}^c(X_\ar,\Z)
\to H_{i+1}^K(X,\Z) \xrightarrow{d_{i-1}} \cdots  $$
for all $X$ (or for all smooth and proper $X$).
In this case, the image of $d_i$ is torsion for all $i$.
\end{theorem}

\proof
If $i\leq d:=\trdeg k/\F_q$, the Beilinson-Lichtenbaum
conjecture implies that $H^i(k,\Z(d))\cong H^i(k_\et,\Z(d))$.
On the other hand, from $H^i(k,\Q(d))=0$ for $i<d$ (known
unconditionally for $d\leq 1$, and for all $d$ by hypothesis in b)), and
the long exact sequence \eqref{alt}, it follows that
$H^i(k_\et,\Z(d))=H^i(k_\ar,\Z(d))$ for $i\leq d$.
Both exact sequences follow by comparing the niveau spectral sequences
\eqref{niveau} and \eqref{niveau1}, because they agree (in the given range)
up to the line defining Kato-homology.
The statement on $d_i$ follows because the canonical map
$CH_0(X,i)_\Q\to H_{i+1}^c(X_\ar,\Q)$ is injective by Theorem \ref{weilI}.
To prove the converse in b), we observe that $H_i^K(X,\Z)$ vanishes for
$i>\dim X$ by definition, hence $CH_0(X,i)=H_{i+1}^c(X_\ar,\Z)$ for
$i\geq \dim X$, and the result follows by tensoring with $\Q$ and
comparing to Theorem \ref{weilI}.
\proofend

\rem
Taking the colimit over smooth schemes with function field $k$,
the same argument shows that the statement in Prop. \ref{eqeqe}a)
implies that for a field $k$ of transcendence degree $d$ over $\F_q$, we have
$$ H^i(k_\ar,\Z(d))\cong \begin{cases}
H^i(k,\Z(d)) &i\leq d;\\
H_d^K(k,\Z)  &i=d+1;\\
0              &i>d+1.
\end{cases}$$

\subsection{The conjecture}
The following is an integral version of Kato's conjecture.

\begin{conjecture}\label{gkato}
If $X$ is smooth, proper and connected, then
$H_i^K(X,\Z)=0$ for $i>0$, and the augmentation map induces
an isomorphism $H_0^K(X,\Z)\cong\Z$.
\end{conjecture}

\begin{proposition}\label{ppplll}
Conjecture \ref{gkato} is equivalent to the conjunction of
Conjecture \ref{katororign} and the statement in Proposition \ref{eqeqe}b).
\end{proposition}

\proof
If Conjecture \ref{gkato} holds, then for smooth and proper $X$, the
groups $H_i^K(X,\Z/m)$  vanish for $i>0$ by \eqref{coecoe}, and then
$\tilde H_i(X,\Q(0))=0$ for $i>0$ and $\tilde H_0(X,\Q(0))=\Q$ by
\eqref{relate}. The converse follows by \eqref{relate}.
\proofend

Note that Propositions \ref{ppplll} and \ref{eqeqe} together with
Theorem \ref{maincft} imply Theorem \ref{intro2}.

Assuming resolution of singularities, Jannsen \cite[Thm.5.9]{jannsen}
defines  weight-homology $H_i^W(X,\Z)$ of $X$ as the homology of the complex
$$ \cdots \to\Z^{\pi_0(X_1)} \to \Z^{\pi_0(X_0)}  ,$$
where
$$ \cdots \to M(X_1)\to  M(X_0)$$
is the weight complex $W(X)$ of Gillet-Soul\'e \cite{gilletsoule},
see the discussion in \cite{ichparshin}.
By definition for $X$ smooth and proper,
$H_0^W(X,\Z)=\Z$, and $H_i^W(X,\Z)=0$ for $i>0$.
If follows from the properties of weight complexes that
$H_*^W(-,\Z)$ has the localization property. We define weight homology
of a field $k$ of finite type over $\F_q$ as
$H_i^W(k,\Z):=\colim H_i^W(U,\Z)$, where $U$ runs through the filtered
system of schemes over $\F_q$ with function field $k$.
By Jannsen \cite[Thm.5.9]{jannsen}, $H_i^W(k,\Z)=0$ for $i\not=\trdeg_{\F_q}k$,
hence the niveau spectral sequence shows that
the weight homology $H_i^W(X,\Z)$ of $X$ is the homology of the complex
\begin{equation}\label{jannsenw}
\cdots \to \oplus_{x\in X_{(s)}}H_s^W(k(x),\Z)
\to \cdots\to \oplus_{x\in X_{(0)}}H_0^W(k(x),\Z).
\end{equation}

\begin{lemma}\label{popob}
For any scheme $X$ over $\F_q$, we have a canonical isomorphism
$$H^K_i(W(X),\Z)\xrightarrow{\sim} H^K_i(X,\Z),$$
compatible with localization triangles.
\end{lemma}

\proof
If $X$ is proper, then $W(X)$ is chain homotopy equivalent to
a hyperenvelope $X_\cdot$ of $X$, and the augmentation
$W(X)\cong X_\cdot\to X$ induces the indicated quasi-isomorphism by the
localization property of Kato homology by the argument of \cite{gillet}.
For an arbitrary scheme $U$, we choose a compactification $X$ and complement
$Z$, and obtain a commutative diagram
$$ \begin{CD}
C^K(W(Z))@>>> C^K(W(X))@>>> C^K(W(U))\\
@V\cong VV @V\cong VV \\
C^K(Z)@>>> C^K(X)@>>> C^K(U).
\end{CD}$$
The rows are distinguished triangles by the localization property
of Kato-homology and property d) of weight complexes, respectively.
\proofend

The augmentation maps $C^K(X_i)\to \Z^{\pi_0(X_i)}$
induce a natural homomorphism
$$ H^K_i(X,\Z)\cong H^K_i(W(X),\Z)\to H^W_i(X,\Z).$$
Taking the colimit over smooth $U$ with function
field $k$, we obtain a map
$H^{d+1}(k_\ar,\Z(d))\cong H_d^K(k,\Z)\to H_d^W(k,\Z)$.

\begin{proposition}\label{lllppp}
Assume resolution of singularities. Then the following statements are
equivalent:
\begin{enumerate}
\item Conjecture \ref{gkato}.
\item The canonical map $H^{d+1}(k_\ar,\Z(d))\to H_d^W(k,\Z)$ is an
isomorphism for all fields $k$ with $\trdeg k/\F_q=d$.
\item The canonical map $H_i^K(X,\Z)\to H_i^W(X,\Z)$ is an isomorphism
for all schemes $X$ and all $i$.
\end{enumerate}
\end{proposition}

\proof
a) $\Rightarrow$ b): We proceed by induction on $d$, the case $d=0$
being trivial. Choose a smooth and proper model $X$ for $k$
and compare the exact sequences \eqref{gkatos} and \eqref{jannsenw}
$$\begin{CD}
H^{d+1}(k_\ar,\Z(d))@>>> \oplus_{x\in X_{(d-1)}}H^d(k(x)_\ar,\Z(d-1))
@>>>\cdots \\
@VVV @|@| \\
H_d^W(k,\Z)@>>>\oplus_{x\in X_{(d-1)}}H_{d-1}^W(k(x),\Z)
@>>>\cdots
\end{CD}$$
The upper sequence is exact by hypothesis, and the lower sequence
is exact by Jannsen's theorem, because $H^W_i(X,\Z)=0$ for $i>0$.

\noindent b) $\Rightarrow$ c) follows by comparing \eqref{gkatos}
and \eqref{jannsenw}, and c) $\Rightarrow$ a) is trivial.
\proofend

\begin{corollary}
Assume resolution of singularities. Then Conjecture $P_0(X)$ holds for
all smooth and proper
$X$ if and only if there is an exact sequence of finitely generated groups
$$\cdots \to CH_0(X,i)\to H_{i+1}^c(X_\ar,\Z)
\to H_{i+1}^W(X,\Z) \to \cdots $$
for all $X$.
\end{corollary}

\proof
If $P_0(X)$ holds for all $X$, then we get the exact sequence
from Theorem \ref{maincft}b), and Prop. \ref{eqeqe}, \ref{ppplll}
and \ref{lllppp}. The finite generation statement follows because weight
homology is finitely generated by construction.
Conversely, the exact sequence gives $CH_0(X,i)\cong H_{i+1}^c(X_\ar,\Z)$
for smooth and proper $X$ and $i\geq 0$, and hence $P_0(X)$ by comparing to
Thm. \ref{weilI}.
\proofend

\section{Curves}

\begin{theorem}
Conjecture \ref{gkato} holds for all smooth and proper curves over $\F_q$.
\end{theorem}

\proof
Consider the following map of short exact sequences
\begin{equation}
\label{curved}\begin{CD}
H^2(k(C)_\ar,\Z(1))@>>>H^1(k(C),\Q(1))@>>> \Br(k(C))\\
@VVV @VVV @VVV \\
\bigoplus_{x\in C_{(0)}}\Z@>>> \bigoplus_{x\in C_{(0)}} \Q @>>>
\bigoplus_{x\in C_{(0)}} \Q/\Z .
\end{CD}\end{equation}
The upper row is \eqref{ikmn}, since $\Br(k\otimes_{\F_q}\bar\F_q)$
vanishes by Tsen's theorem.
The snake Lemma gives a short exact sequence of kernels and cokernels
$$ 0\to H_1^K(C,\Z)\to  K_1(C)_\Q\to \Br(C)\to H_0^K(C,\Z)\to
K_0(C)_\Q  \to \Q/\Z\to 0,$$
and the claim follows because $K_1(C)_\Q=0=\Br(C)$.
\proofend

The above theorem can be derived
without using the vanishing of the Brauer group by analyzing
the Frobenius coinvariants of the localization sequence
$$ 0\to \Gamma(\bar C,{\mathcal O}_{\bar C})^\times
\to (k(C)\otimes_{\F_q}\bar\F_q)^\times \to \oplus_{x\in \bar C_{(0)}}\Z
\to \Pic \bar C\to 0.$$
It is amusing to observe that
$$ 0\to H^2(k(C)_\ar,\Z(1))\to \bigoplus_{x\in C_{(0)}}\Z\to \Z\to 0$$
is an integral version of the classical short exact sequence
$$ 0\to \Br(k(C))\to \bigoplus_{x\in C_{(0)}} \Q/\Z \to \Q/\Z\to0.$$
Indeed, the upper row of \eqref{curved}
shows that $H^2(k(C)_\ar,\Z(1))\otimes\Q/\Z\cong \Br(k(C))$.

\begin{theorem}\label{curve}
Let $C$ be any curve. Then $H_0^c(C_\ar,\Z)\cong H_0^K(C,\Z)$,
there is a short exact sequence
$$ 0\to CH_0(C)\to H_1^c(C_\ar,\Z)\to H_1^K(C,\Z)\to 0,$$
an isomorphism of finitely
generated groups $CH_0(C,1)\cong H_2^c(C_\ar,\Z)$,
and $CH_0(C,i)=H_{i+1}^c(C_\ar,\Z)=0$ for $i>1$.
\end{theorem}

The Pointrjagin dual of the short exact sequence is
\cite[Prop. 1]{katosaito}.

\proof
The proof of Thm. \ref{maincft} works if one
restricts oneself to schemes of dimension at most $1$, in which
case the Beilinson-Lichtenbaum conjecture and $P_0(X)$ is known.
Similarly, finite generation follows by Prop. \ref{4.1}
and \ref{5.3}.
\proofend

For a curve $C$, integral Kato homology can be calculated
combinatorically. We first assume that $C$ is connected and proper.
Let $C'=\coprod_{i\in I} C_i'$
be the decomposition of the normalization $C'$ of $C$ into
irreducible components, $S$ the set of singular points of $C$,
and $S'=S\times_CC'$.
Consider the dual graph $\Gamma$ of $C$ \cite{roberts}.
It is a bipartite graph,
with vertices the set $S$ of singular points of $C$ on the one hand,
and the set of irreducible components $I$ of $C'$ on the other hand. For each
point $s\in S'$, there is an edge connecting the image of $s$ in $S$
with the $i\in I$ such that $s\in C_i$.
Comparing the exact sequence calculating the homology of $\Gamma$
to the sequence
$$ 0\to H_1^K(C,\Z) \to H_0^K(S',\Z)\to H_0^K(S,\Z)\oplus H_0^K(C',\Z)
\to H_0^K(C,\Z)\to 0,$$
we get

\begin{proposition}\label{curves}
If $C$ is a proper curve, then $H_i^K(C,\Z)\cong H_i(\Gamma,\Z)$.
In particular,
$H_0^K(C,\Z)=\Z^{\pi_0(C)}$,
and $H_1^K(C,\Z)= \Z^{\pi_0(S')-\pi_0(S)-\pi_0(C')+\pi_0(C)}$.
\end{proposition}

If $C$ has only ordinary
double points as singularities, then the homology of $\Gamma$
agrees with the homology of the graph $\Gamma'$
used in \cite{ams}. Indeed, the graph
$\Gamma'$ is gotten from $\Gamma$ by erasing all vertices corresponding
to points in $S$. Since over every point of $s\in S$ there are exactly
two point of $x,y\in S'$, the two edges corresponding to $x$ and $y$
combine to give an edge labeled $s$ in $\Gamma'$. Thus $\Gamma$ and
$\Gamma'$ are homeomorphic as topological spaces.

If $C$ is not proper, let $\bar C$ be a proper curve containing $C$.
Then the exact sequence
$$ 0\to H_1^K(\bar C,\Z)\to H_1^K(C,\Z)\to \bigoplus_{p\in \bar C-C}\Z
\to H_0^K(\bar C,\Z)\to H_0^K(C ,\Z)\to 0$$
shows that $H_0^K(C,\Z)=0$ and
$H_1^K(C,\Z)=H_1^K(\bar C,\Z)\oplus \Z^{\pi_0(\bar C-C)-1}$.

\section{Class field theory}
In this section, we assume that $X$ is proper and connected over $\F_q$.

\begin{proposition}\label{ghy}
We have
$$ H_1^c(X_\ar,\Z)^\wedge \cong \pi_1^{ab}(X).$$
In particular, if $H_1^c(X_\ar,\Z)$ is finitely generated,
then there is a short exact sequence
$$ 0\to (\CH_0(\bar X,1)_G)^\wedge \to \pi_1^{ab}(X) \to
(\CH_0(\bar X)^G)^\wedge \to  0.$$
\end{proposition}

\proof
By Prop. \ref{degreenull}, $H_0^c(X_\ar,\Z)$ is
torsion free, hence by Thm. \ref{weilI}
$$H_1^c(X_\ar,\Z)/m\cong H_1^c(X_\ar,\Z/m)\cong H_0(X_\et,\Z/m).$$
By \cite{ichdual}, the latter is isomorphic to
$H^1(X_\et,\Z/m)^*=\pi_1^{ab}(X)/m$.
The Proposition follows by taking inverse limits.
\proofend

The results of this section could be formulated independently of
Conjecture $P_0(X)$ using the dual of arithmetic cohomology of
\cite{ichweilII}:

\begin{proposition}
We have a natural inclusion
$$H^{-1}R\Hom_\Ab(R\Gamma(X_\ar,\Z),\Z) \to \pi_1^{ab}(X)$$
which induces an isomorphism on completions.
There are short exact sequences
$$ 0\to (H^2(X_\ar,\Z)_\tor)^* \to \pi_1(X)^{ab}\to
\Hom(H^1(X_\ar,\Z),\Z)^\wedge\to 0.$$
For normal $X$, under resolution of singularities,
$(H^2(X_\ar,\Z)_\tor)^*$ is isomorphic to the geometric part of the
abelianized fundamental group $\pi_1^{ab}(X)^0$.
\end{proposition}

\proof
Consider the reduction mod $m$ map
\begin{multline*}
H^{-1}R\Hom_\Ab(R\Gamma(X_\ar,\Z),\Z)/m \to
H^{-1}(R\Hom_\Ab(R\Gamma(X_\et,\Z),\Z)\otimes^L\Z/m)\\
\to H^{-1}(R\Hom_\Ab(R\Gamma(X_\et,\Z/m),\Z/m))
\cong \Hom(H^1(X_\ar,\Z/m),\Q/\Z).
\end{multline*}
Taking the limit over $m$, we obtain a map
\begin{multline*}
$$H^{-1}R\Hom_\Ab(R\Gamma(X_\ar,\Z),\Z)\to
\lim\Hom(H^1(X_\ar,\Z/m),\Q/\Z) \\
=\Hom(H^1(X_\ar,\Q/\Z),\Q/\Z)
=\Hom(H^1(X_\et,\Q/\Z),\Q/\Z)=\pi_1(X)^{ab}.
\end{multline*}
Since the first term is finitely generated, it injects
into its completion.
The short exact sequence is the Leray spectral sequence
for $R\Hom_\Ab(-,\Z)$ (using that
$\Ext^1(-,\Z)\cong (A_\tor)^*$ for finitely generated $A$).
For normal $X$, $H^1(\bar X_\eh,\Z)=0$ under resolution of
singularities \cite{ichweilII},
hence $H^1(X_\ar,\Z)\cong H^0(\bar X_\eh,\Z)_G\cong \Z$.
\proofend

Let $X$ be smooth, projective and geometrically connected.
The isomorphism of class field theory of Bloch and
Kato-Saito \cite{katosaito}, factors through arithmetic homology
$$\theta_X:
CH_0(X)_\tor\to H_1^c(X_\ar,\Z)_\tor\cong \pi_1^{ab}(X)_\tor.$$
Let $T$ be the finite group $\colim\Hom_{GS/k}(\mu_m,\NS_X)$,
where $\Hom_{GS/k}(-,-)$ denotes the group of homomorphisms of groups
schemes over $k$, and recall that conjecturally $CH_0(X,1)_\Q=0.$

\begin{proposition}\label{smoothcase}
There are short exact sequences
\begin{align*}
0\to \CH_0(\bar X,1)_G\to &H_1^c(X_\ar,\Z)\to \Alb(X)\oplus \Z\to 0,\\
0\to CH_0(X)\to &H_1^c(X_\ar,\Z)\to CH_0(X,1)_\Q \to 0,
\end{align*}
and
$$0\to T^* \to CH_0(\bar X,1)_G \to CH_0(X,1)_\Q\to 0.$$
\end{proposition}

\proof
The first exact sequence is \eqref{lerayses}, together with
$\CH_0(\bar X)^G\cong (\Alb(\bar X)\oplus\Z)^G=\Alb(X)\oplus \Z$.
The next sequence follows from Cor. \ref{known} and
Thm. \ref{maincft}. Indeed, $H_2^K(X,\Z)$ is uniquely
divisible, hence its image in the finitely generated group
$CH_0(X)$ must be trivial.
The last sequence follows by considering the following diagram
(the upper row is \cite[Prop. 9(2)]{katosaito})
$$\begin{CD}
0@>>>  T^* @>>>  CH_0(X) @>>> CH_0(\bar X)^G@>>>  0\\
@III @VVV @VVV @| \\
0@>>> CH_0(\bar X,1)_G@>>> H_1^c(X_\ar,\Z)@>>> CH_0(\bar X)^G@>>> 0.
\end{CD}$$
\proofend

Consider a cartesian diagram
$$\begin{CD}
Z'@>>> X'\\
@VVV @VpVV\\
Z@>i>> X;
\end{CD}$$
such that  $i$ is a closed embedding, $p$ is proper, and $p$ induces an
isomorphism $X'-Z'\cong X-Z$.

\begin{proposition}\label{general}
If $X'$ and $Z$ are smooth, then  there is an exact sequence
$$ \ker \big(H_1^K(Z',\Z)\to H_1^K(X',\Z)\oplus H_1^K(Z,\Z)\big)  \to \CH_0(X)
\stackrel{\theta_X}{\longrightarrow} H_1^c(X_\ar,\Z).$$
The first term agrees with $H_1^K(Z',\Z)$ if
$CH_0(X',1)_\Q\cong CH_0(Z,1)_\Q=0$.
\end{proposition}

\proof
Consider the commutative diagram with
exact rows and exact columns given by Theorem \ref{maincft},
$$\begin{CD}
\CH_0(Z')@>>> \CH_0(X')\oplus \CH_0(Z) @>>> \CH_0(X) \to 0\\
@V\theta_{Z'}VV @V\theta_{X'}V\theta_ZV @V\theta_X VV  \\
H_1^c(Z'_\ar,\Z)@>>>  H_1^c(X'_\ar,\Z)\oplus H_1^c(Z_\ar,\Z)@>>>
H_1^c(X_\ar,\Z)\\
@VVV @VVV @VVV\\
H_1^K(Z',\Z)@>>>  H_1^K(X',\Z)\oplus H_1^K(Z,\Z)@>>> H_1^K(X,\Z)\\
@VVV @VVV @VVV \\
0@EEE 0@EEE 0.
\end{CD}$$
Using the injectivity of $\theta_{X'}$ and $\theta_Z$, the first
statement follows by a diagram chase.
The last claim follows from Corollary \ref{known}.
\proofend

Theorem \ref{maincft} and Proposition \ref{general} can be applied
to study the kernel of the reciprocity map.
For example, we can recover some results of Asakura, Matsumi and Sato
\cite{ams}. Let $X$ be a normal surface with one singular point $P$.
Let $S$ be a resolution of singularities $X$ which is an isomorphism
away from $P$, and let $D=P\times_X S$ be the exceptional fiber.

\begin{proposition}
Let $\Gamma$ be the dual graph of $D$, and
assume that $CH_0(S,1)_\Q=0$. Then 
there are short exact sequences
$$ H_1(\Gamma,\Z)\stackrel{\delta_X}{\longrightarrow}
CH_0(X)\stackrel{\theta_X}{\longrightarrow}
H_1^c(X_\ar,\Z),$$
and
\begin{multline*}
0\to \ker\big(CH_0(S)\to CH_0(X)\big)\to
\ker\big(H_1^c(S_\ar,\Z)\to H_1^c(X_\ar,\Z)\big)\\
 \to \ker \theta_X
\to \coker\big(CH_0(D)\stackrel{\deg}{\to} \Z\big)\to  0.
\end{multline*}
\end{proposition}

\proof
The first statement follows from Props. \ref{curves}
and \ref{general}. To prove the second statement, consider the commutative
diagram
\begin{equation}\label{thedg}
\begin{CD}
@EEE 0@>>> CH_0(S)@>\theta_S>\sim> H_1^c(S_\ar,\Z)@>>> 0\\
@III @VVV @VVV@VVV \\
0@>>> \ker \theta_X @>>>CH_0(X)@>\theta_X>>H_1^c(X_\ar,\Z)@>>> 0.
\end{CD}
\end{equation}
The upper row is exact by Prop. \ref{smoothcase}. Since
$\coker\big( CH_0(S)\to CH_0(X)\big) =\coker\big( CH_0(D)\to CH_0(P)\big)$,
it suffices to show that $H_1^c(S_\ar,\Z)\to H_1^c(X_\ar,\Z)$ is surjective.
Since $X$ is normal, $D$ is connected, hence
$H_0^c(D_\ar,\Z)\cong H_0^c(P_\ar,\Z)$, and we obtain an exact sequence
$$ H_1^c(D_\ar,\Z)\to H_1^c(P_\ar,\Z)\oplus H_1^c(S_\ar,\Z)
\to H_1^c(X_\ar,\Z) \to 0.$$
Thus it suffices to show that $H_1^c(D_\ar,\Z)\to H_1^c(P_\ar,\Z)$ is surjective.
Since $H_1^c(P_\ar,\Z)\cong \Z$, we can show surjectivity
modulo $m$ for every $m$. By torsion-freeness of $H_0^c(D_\ar,\Z)$,
we obtain a diagram
$$\begin{CD}
H_1^c(D_\ar,\Z)/m @= H_1^c(D_\ar,\Z/m) @= H^1(D_\et,\Z/m)^*\\
@VVV @VVV @VVV \\
H_1^c(P_\ar,\Z)/m @= H_1^c(P_\ar,\Z/m) @= H^1(P_\et,\Z/m)^*
\end{CD}$$
and can show injectivity of
$H^1(P_\et,\Z/m)\to H^1(D_\et,\Z/m)$. By the Hochschild-Serre
spectral sequence, this  follows from injectivity of
$H^1(\hat G,H^0(\bar P_\et,\Z/m))\to H^1(\hat G,H^0(\bar D_\et,\Z/m))$.
But the latter map is an isomorphism because $\bar D$ is connected
by Zariski's main theorem.
\proofend

According to \cite[Lemma 1.2, Prop. B.6]{ams},
$\coker \big(CH_0(D)\to \Z\big)\cong \Z/m$,
where $m$ is the greatest common
divisor of the degree of the field extension
$\Gamma(C_i,{\mathcal O}_{C_i})/\Gamma(D,{\mathcal O}_D)$.
Here $C_i$ runs through the irreducible components of $D$.

\end{document}